\documentclass[12pt]{article}
\usepackage{amsmath,amssymb,amsthm,amsfonts,amsxtra,yhmath}
\newtheorem{thm}{Theorem} 
\newtheorem{prop}[thm]{Proposition}
\newtheorem{lemma}[thm]{Lemma}
\newtheorem{cor}[thm]{Corollary}
\theoremstyle{definition}

\theoremstyle{remark}
\newtheorem{rem}[thm]{Remark}
\newtheorem{example}[thm]{Example}

\newcommand{\set}[2]{ \left\{\,#1\,;\,#2\,\right\} }
\newcommand{\cpling}[2]{ \langle #1,#2 \rangle}
\newcommand{\transp}[1]{ {}^{\mathsf{t}} #1}

\newcommand{\changed}[1]{#1}

\newcommand{\nc}{\newcommand}

\nc{\sh}{y}
\nc{\g}{\mathfrak{g}}
\nc{\ep}{\varepsilon}
\nc{\iu}{i}
\nc{\integer}{\mathbb{Z}}
\nc{\real}{\mathbb{R}}
\nc{\complex}{\mathbb{C}}
\nc{\SH}{\mathcal{SH}}

\nc{\shdual}{z} 
\nc{\h}{\mathfrak{h}}
\nc{\n}{\mathfrak{n}}
\nc{\us}{\underline{s}}
\nc{\Inv}{\mathcal{I}}
\nc{\ML}{\mathcal{ML}}

\begin{document}
\title{Some properties of associated spaces with sub-Hankel determinants}
\author{%
Hideyuki Ishi%
\thanks{The second authors are partially supported by the grant in aid of scientific research of JSPS 24540049,}
\footnote{Graduate School of Mathematics,
Nagoya University, Furo-cho, Nagoya,  464-8602,
Japan. \newline e-mail: \texttt{hideyuki@math.nagoya-u.ac.jp }}
\quad and \quad Takeyoshi Kogiso${}^*$\!\!
\footnote{Department of Mathematics, Josai University, 1-1 Keyakidai, Sakado, Saitama, 350-0295, Japan.\newline e-mail: \texttt{kogiso@math.josai.ac.jp}}%
}
\date{}

\maketitle

\begin{abstract} 
 In this note, we show that the space associated with sub-Hankel determinant is
 a non-reductive, regular prehomogeneous vector space, and we give the
 multiplicative Legendre transforms of sub-Hankel determinants. Moreover we 
observe 
  {certain relations}
between b-functions of polarization  of  PV-polynomials and $b$-functions of sub-Hankel determinants,
 and give some formulas 
 about sub-Hankel determinants whose components are orthogonal ponlynomials.
\end{abstract}

\section*{Introduction}
In the paper \cite{CRS}, the notion of sub-Hankel determinant has been introduced 
as an interesting example of homaloidal polynomial.
According to \cite{CRS}, the sub-Hankel matrix $M^{(r)}$ of size $r$ is a matrix whose $(i,j)$-component depends only on $i+j$,
 and equals $0$ if $i + j \ge r+2$.
For sub-Hankel determinants, 
 we shall present in this paper the following results and conjectures:

i) The space $\SH(r)$ associated with the sub-Hankel determinant $\det M ^ {(r)} $ is a regular prehomogeneous vector space,
 and $P_1 = \det M^{(r)}$ is one of the two basic relative 
invariants $P_i\,\,\,(i=1,2)$ of the space.

ii) We determine two basic relative invariants $Q_i\,\,(i=1,2)$ of the dual prehomogenous vector space 
$\SH(r)^*$
of the space $\SH(r)$,
and compute the multiplicative Legendre transforms of 
generic $P_1^{\alpha_1} P_2^{\alpha_2}$ and $Q_1^{\beta_1}Q_2^{\beta_2}$
respectively. 

iii) By comparing the corresponding rational characters of $P_1, P_2, Q_1,Q_2$, we pose a conjecture about the form 
of $b$-function of   {$\det M^{(r)}$}.

iv) We introduce some formulas obtained by substituting some kinds of orthogonal polynomials.



\section{Structures of associated spaces to sub-Hankel determinants }
Let $\SH(r)$ be the vector space of sub-Hankel matrices of size $r$.
Namely,
 $\SH(r)$ is the set of matrices $\sh$ of the following form:
\begin{equation} \label{eqn:subHankel_y}
\sh = \begin{pmatrix}
 y_1 &        y_2 &   y_3 & \cdots & y_{r-1} &  y_r \\
 y_2 &        y_3 &   y_4 & \cdots  & y_r    & y_{r+1} \\
 y_3 &        y_4 &   y_5    &    \cdots  & y_{r+1} & 0 \\
 \vdots & \vdots  & \adots & \adots &    0  & \vdots\\
 y_{r-1}  &     y_r& y_{r+1} & \adots & \vdots& 0 \\
 y_r     & y_{r+1} &     0 & \cdots &     0 & 0  
\end{pmatrix}.
\end{equation}
For $\ell=1, \dots, r+1$,
 let $Y_{\ell} \in \mathrm{Mat}(r, \real)$ be the matrix whose $(i,j)$-component is $\delta_{i+j,\, \ell+1}$.
Then $Y_1, \dots, Y_{r+1}$ forms a basis of $\SH(r)$. 
Indeed,
 the element $\sh \in \SH(r)$ in (\ref{eqn:subHankel_y}) equals $\sum_{\ell=1}^{r+1} y_{\ell} Y_{\ell}$.
For $k=0, \dots, r-1$,
 we write $\mathrm{diag}^{(k)}(a_1, \dots, a_{r-k})$ for the $r \times r$ upper triangular matrix
 whose $(i,j)$-component is $a_i \delta_{i,\, j-k}$.
We define 
$$ T_k := \mathrm{diag}^{(k)}(r-k, r-k-1, \dots, 1). $$


\begin{lemma} \label{lemma:Tk-Yl}
{\rm (i)} For $k, k' =0, 1,\dots, r-1$, 
 one has
$$
 [T_k, T_{k'}] = \begin{cases} (k'-k) T_{k+k'} & (\mbox{if }0 \le k+k' \le r-1), \\ 0 & (\mbox{otherwise}). \end{cases}
$$ 
{\rm (ii)} For $k=0, \dots, r-1$ and $\ell=1, \dots, r+1$, 
one has
$$
 T_k Y_{\ell} + Y_{\ell} \transp{T_k} 
 = \begin{cases} (2r+1 - k - \ell) Y_{\ell - k} & (\mbox{if }\ell-k \ge 1), \\ 0 & (\mbox{otherwise}). \end{cases}
$$
\end{lemma}

 
Let us define
\begin{align*} 
 H_1 &:= \frac{1}{r} T_0 - \frac{1}{2} I_r 
 = \frac{1}{2} I_r - \mathrm{diag}(0, \frac{1}{r}, \dots, \frac{r-1}{r}).\\
 H_2 &:= I_r - \frac{1}{r} T_0 
 = \mathrm{diag}(0, \frac{1}{r}, \dots, \frac{r-1}{r}).
\end{align*}
We see from Lemma~\ref{lemma:Tk-Yl} that
\begin{equation} \label{eqn:HTk}
[H_1, T_k] = \frac{k}{r} T_k, \quad
[H_2, T_k] =-\frac{k}{r} T_k \qquad (k=1, \dots, r-1),
\end{equation}
and that
\begin{equation} \label{eqn:HYl}
H_1 Y_{\ell} + Y_{\ell} \transp{H_1} = \frac{r+1 - \ell}{r} Y_{\ell}, \quad
H_2 Y_{\ell} + Y_{\ell} \transp{H_2} = \frac{\ell -1}{r} Y_{\ell} \qquad (\ell = 1, \dots, r+1).
\end{equation} 
Let $\g \subset \mathrm{Mat}(r, \real)$ be the vector space spanned by $H_1, H_2, T_1, \dots, T_{r-1}$.
We see from Lemma~\ref{lemma:Tk-Yl} (i) and (\ref{eqn:HTk}) that $\g$ forms a solvable Lie algebra. 
Let $G\subset GL(r, \real)$
 be the Lie group $\exp \g$ corresponding to $\g$.
Thanks to Lemma~\ref{lemma:Tk-Yl} (ii) and (\ref{eqn:HYl}),
 we can define a representation $\rho$ of the group $G$ on the space $\SH(r)$ by
$$
 \rho(\exp T) \sh := (\exp T) \sh\, \transp{(\exp T)} \qquad (T \in \g)
$$
 so that its infinitesimal representation is given by
$$
d\rho(T) \sh := T \sh + \sh \, \transp{T} \qquad (T \in \g).
$$

\section{Multiplicative Legendre transforms of sub-Hankel determinants}

For $\us = (s_1, s_2) \in \complex^2$,
 let $\nu_{\us} : G \to \complex^{\times}$ be the one-dimensional representation
 of the group $G$ given by
$$
\nu_{\us}(\exp T) := e^{s_1 h_1 + s_2 h_2} \quad (T= \sum_{i=1}^2 h_i H_i + \sum_{k=1}^{r-1} t_k T_k \in \g).
$$
For $a = \exp(h_1 H_1 + h_2 H_2) \in G\,\,\,(h_1, h_2 \in \real)$, 
 we have
\begin{equation} \label{eqn:expH}
 a = \mathrm{diag}(a_1^{1/2},\, a_1^{1/2 - 1/r},\, a_1^{1/2 - 2/r}, \dots, a_1^{1/2 -(r- 1)/r})  
          \cdot \mathrm{diag}(1,\, a_2^{1/r}, a_2^{2/r}, \dots, a_2^{(r-1)/r}),
\end{equation} 
 where $a_i := e^{h_i}$ for $i=1,2$.
Then we have 
 $$\nu_{\us}(a) = a_1^{s_1} a_2^{s_2}.$$

In view of the action $\rho$,
 we see that the right-bottom principal minors of $\sh$
 as well as the determinant $\det \sh$
 are $\rho(G)$-relatively invariant.
Note that the right-bottom principal minor of degree $r-1$ equals $(-1)^{[\frac{r-1}{2}]}(y_{r+1})^{r-1}$,
 \changed{where $[\ \cdot \ ]$ is the Gauss symbol.}
We define
\begin{equation} \label{eqn:def_of_Pi}
 P_1(\sh) := \det \sh, \quad P_2(\sh) := y_{k+1} \qquad (\sh \in \SH(r)),
\end{equation}
 which are $\rho(G)$-relatively invariant polynomial functions on $\SH(r)$.
For the element $a$ in (\ref{eqn:expH}), we have
$\det(\rho(a)\sh) = (\det a)^2 \det \sh = a_1 a_2^{r-1} \det \sh$.
Thus
\begin{equation} \label{eqn:P1}
P_1(\rho(g)\sh) = \nu_{(1,r-1)}(g) P_1(\sh) \quad (g \in G,\,\sh \in \SH(r)).
\end{equation}
On the other hand, we have $P_2(\rho(a)\sh) = a_2 y_{r+1} = a_2 P_2(\sh)$, so that 
\begin{equation} \label{eqn:P2}
P_2(\rho(g)\sh) = \nu_{(0,1)}(g) P_2(\sh) \quad (g \in G,\,\sh \in \SH(r)).
\end{equation}


\begin{prop} \label{prop:y-orbit}
{\rm (i)}If $P_1(\sh) \ne 0$ and $P_2(\sh) \ne 0$, then the orbit $\rho(G)\sh$ in $\SH(r)$ is described as
$$
 \rho(G)\sh = \set{\sh' \in \SH(r)}{P_1(\sh)P_1(\sh')>0 \mbox{ and }P_2(\sh)P_2(\sh')>0}.
$$
{\rm (ii)} The triple $(G, \rho, \SH(r))$ is a prehomogeneous vector space whose singular set is
$\Sigma := \set{\sh \in \SH(r)}{P_1(\sh)  P_2(\sh) = 0}$.
\end{prop}


For $\us =(s_1, s_2) \in \complex^2$,
 we put
$$ \phi_{\us}(\sh) 
 := \Bigl\{ \frac{\det \sh}{(y_{r+1})^{r-1}} \Bigr\}^{s_1} (y_{r+1})^{s_2} 
 = P_1(\sh)^{s_1} P_2(\sh)^{s_2 - (r-1) s_1}.
$$
Then we see from (\ref{eqn:P1}) and (\ref{eqn:P2}) that
\begin{equation} \label{eqn:phi_s}
 \phi_{\us}(\rho(g)\sh) = \nu_{\us}(g) \phi_{\us}(\sh) \qquad (g \in G,\,\sh \in \SH(r)).
\end{equation}

Let $\SH(r)^*$ be the dual vector space of $\SH(r)$.
The contragredient representation $\rho^* : G \to GL(\SH(r)^*)$ is defined by
$$
 \cpling{\sh}{\rho^*(g)\shdual} := \cpling{\rho(g)^{-1}\sh}{\shdual} \qquad (g \in G,\,\sh \in \SH(r),\shdual \in \SH(r)^*).
$$ 
Let $(Y_1^*, \dots, Y_{r+1}^*)$ denote the basis of $\SH(r)^*$
 dual to the basis $(Y_1, \dots, Y_{r+1})$ of $\SH(r)$.
The matrix expression of $\rho^*(g)$
 with respect to $(Y_1^*, \dots, Y_{r+1}^*)$
 is a lower triangular matrix.
Therefore,
 if we define
$$ Q_2 (\shdual) := \shdual_1 \qquad (\shdual = \sum_{\ell=1}^{r+1} \shdual_{\ell} Y_{\ell}^* \in \SH(r)^*,\, \shdual_1, \dots, \shdual_{r+1} \in \real), $$
 then $Q_2 (\rho^*(a)\shdual) = a_1 \shdual_1 = a_1 Q_2 (\shdual)$
 for $a$ in (\ref{eqn:expH}),
 so that
\begin{equation} \label{eqn:Q1}
   {Q_2}(\rho^*(g)\shdual) = \nu_{(1,0)}(g)   {Q_2}(\shdual) \qquad (g \in G,\, \shdual \in \SH(r)^*).
\end{equation}

Let us define $\g' := \mathrm{span}\langle H_2, T_1, \dots, T_{r-1}\rangle$.
Then $\g'$ is an ideal of $\g$, and 
  {we see that the space
$d\rho^*(\g') (\SH(r)^*)$ equals $W := \mathrm{span}\langle Y^*_2, \dots, Y^*_{r+1} \rangle$.}
We introduce a linear map $R: \SH(r)^* \to \mathrm{Hom}(\g', W)$ defined by
$$
 R(\shdual) T := d\rho^*(T) \shdual \in W \qquad (\shdual \in   {\SH(r)^*},\, T \in \g').
$$
Then we have for $g \in G$
\begin{equation} \label{eqn:equiv_R}
 R(\rho^*(g) \shdual) = \rho^*(g)|_{W} \circ R(\shdual) \circ \mathrm{Ad}(g^{-1})|_{\g'}
\end{equation} 
because 
$$
 R(\rho^*(g) \shdual)X 
= d\rho^*(X) \circ \rho^*(g) \shdual 
= \rho^*(g) \circ d\rho^*(\mathrm{Ad}(g^{-1})X) \shdual
= \rho^*(g) \circ R(\shdual) \circ \mathrm{Ad}(g^{-1}) X 
$$
 for $X \in \g'$.
We   {put}
 $Q_1 (\shdual) := \changed{\frac{(-1)^{r+1}\, 2^{\ell_2(r)}}{r!}}  \det R(\shdual)\quad (\shdual \in \SH(r)^*)$, 
 where the determinant is defined with the bases 
 $(T_{r-1}, \dots, T_1, r H_2)$ of $\g'$ and $(-Y^*_2, \dots, -Y^*_{r+1} )$ of $W$,
 \changed{and $2^{\ell_2 (r) } $ means that maximal 2 power integer factors of $r!$,  
 that is,
 $\ell_2 (r)=\sum_{k \geq 1} \left[ \frac{r}{2^k} \right]$. 
The coefficient is taken for a convenience of normalization.}
We can compute
 $\det \rho^*(g)|_{W} = \nu_{(-\frac{r-1}{2},\,-\frac{r+1}{2})}(g)$
 and
 $\det \mathrm{Ad}(g)|_{\g'} = \nu_{(\frac{r-1}{2},\, - \frac{r-1}{2}) }(g)$
 for $g \in G$,
 which together with (\ref{eqn:equiv_R}) tells us that
\begin{equation} \label{eqn:Q2}
 Q_1 (\rho^*(g)\shdual) = \nu_{(-r+1,-1)}(g) Q_1 (\shdual) \qquad (g \in G,\,\shdual \in \SH(r)^*).
\end{equation}
For $\us =(s_1, s_2) \in \complex^2$,
 we put
$$ \psi_{\us}(\shdual) 
 := (\shdual_1)^{s_1} \Bigl\{ \frac{\det R(\shdual)}{(\shdual_r)^{r-1}} \Bigr\}^{s_2} 
 = Q_2 (\shdual)^{s_1 - (r-1) s_2}  Q_1 (\shdual)^{ s_2} \qquad (\shdual \in \SH(r)^*).
$$
Then we see from (\ref{eqn:Q1}) and (\ref{eqn:Q2}) that
\begin{equation} \label{eqn:psi_s}
 \psi_{\us}(\rho^*(g)\shdual) = \nu_{\us}(g) \psi_{\us}(\shdual) \qquad (g \in G,\,\shdual \in \SH(r)^*).
\end{equation}
Similarly to Proposition \ref{prop:y-orbit}, we have


\begin{prop} \label{prop:z-orbit}
{\rm (i)}If $Q_1(\shdual) \ne 0$ and $Q_2(\shdual) \ne 0$, then the orbit $\rho^*(G)\shdual$ in $\SH(r)^*$ is described as
$$
 \rho^*(G)\shdual = \set{\shdual' \in \SH(r)^*}{Q_1(\shdual)Q_1(\shdual')>0 \mbox{ and }Q_2(\shdual)Q_2(\shdual')>0}.
$$
{\rm (ii)} The triple $(G, \rho^*, \SH(r)^*)$ is a prehomogeneous vector space whose singular set is
$\Sigma^* := \set{\shdual \in \SH(r)^*}{Q_1(\shdual)  Q_2(\shdual) = 0}$.
\end{prop}


For $\us \in \complex^2$ and $\sh \in \SH(r) \setminus \Sigma$,
 we write $\Inv_{\us}(\sh)$ for $\mathrm{grad}\,\log \phi_{\us}(\sh) \in \SH(r)^*$.
Namely,
 we define
$$
 \cpling{v}{\Inv_{\us}(\sh)} := \Bigl(\frac{d}{dt}\Bigr)_{t=0} \log \phi_{\us}(y + t v) \qquad (v \in \SH(r)).
$$
Thanks to the relative invariance (\ref{eqn:phi_s}) of $\phi_{\us}$,
 we have
\begin{equation} \label{eqn:equiv_Inv}
 \Inv_{\us}(\rho(g) \sh) = \rho^*(g) \circ \Inv_{\us}(\sh) \qquad (g \in G).
\end{equation}
On the other hand,
   {we see from some computation that}
\begin{equation} \label{eqn:inv-ab}
 \Inv_{\us}(  {y_1} Y_1 +   {y_{r+1}} Y_{r+1}) = \frac{s_1}{  {y_1}} Y^*_1 + \frac{s_2}{  {y_{r+1}}} Y^*_2
\end{equation}
 for $  {y_1,\, y_{r+1}} \in \real \setminus \{0\}$.
By these observations, we conclude the following.


\begin{prop} \label{prop:bijective_Inv}
The map $\Inv_{\us}$ gives a bijection
 from $\SH(r) \setminus \Sigma$
 onto $\SH(r)^* \setminus \Sigma^*$
 if and only if $s_1 \ne 0$ and $s_2 \ne 0$.
\end{prop}

 
Thanks to Proposition~\ref{prop:bijective_Inv},
 we can define the multiplicative Legendre transform
 $\ML(\phi_{\us})(\shdual):= (\phi_{\us} \circ \Inv_{\us}^{-1}(\shdual))^{-1}$
 as a function on $\SH(r)^* \setminus \Sigma^*$.
For $\shdual \in \SH(r)^* \setminus \Sigma^*$ and $g \in G$,
 we see from (\ref{eqn:equiv_Inv}) and (\ref{eqn:phi_s}) that
\begin{align*}
 \ML(\phi_{\us})(\rho^*(g)\shdual)^{-1} 
&= \phi_{\us} \circ \Inv_{\us}^{-1} \circ \rho^*(g) \shdual
= \phi_{\us} \circ \rho(g) \circ \Inv_{\us}^{-1}(\shdual)
= \nu_{\us}(g) \,\phi_{\us} \circ \Inv_{\us}^{-1}(\shdual)\\
&= \nu_{\us}(g)  \ML(\phi_{\us})(\shdual)^{-1}. 
\end{align*}
Comparing the above with (\ref{eqn:psi_s}),
 we arrive at the following.


\begin{thm} \label{thm:ML_of_phi_s}
If $s_1 \ne 0$ and $s_2 \ne 0$, then the multiplicative Legendre transform $\ML(\phi_{\us})$ is equal to $\psi_{\us}$
 up to constant multiplication.
\end{thm}


\begin{example}
For a sub-Hankel matrix of size 4, that is, 
$y=M^{(4)}=  {\left( \begin{array}{cccc} 
y_{{1}}&y_{{2}}&y_{{3}}&y_{{4}}\\ 
y_{{2}}&y_{{3}}&y_{{4}}&y_{{5}}\\ 
y_{{3}}&y_{{4}}&y_{{5}}&0\\  
y_4 & y_{{5}} & 0 & 0 
\end{array} \right)} $, 
the polynomials
$P_1(y)=\det M^{(4)}$ corresponding to the character 
$\nu_{  {(1, 3)}}$   {and} $P_2 (y)=y_5$ corresponding to the character 
$\nu_{  {(0,1)}}$ are   {the} basic relative invariants of   {the} prehomogeneous vector space 
$( \mathfrak{g}, d \rho , \mathcal{SH}(4))$. 
Then we have
\begin{align*}
   {Q_1(z) =} 
 & \changed{(-1)^5 \det R(z)/3} 
= \changed{-\frac{1}{3}} \left| \begin{array}{cccc} 
 0&0&6\,z_{{1}}& z_{{2}}\\ 
0&4\,z_{{1}}&5\,z_{{2}}& 2\,z_{{3}}
\\ \noalign{\medskip}2\,z_{{1}}&3\,z_{{2}}&4\,z_{{3}}& 3\,z_{{4}}\\
 z_{{2}}&2\,z_{{3}}&3\,z_{{4}}&4\,z_{{5}}
 \end{array} \right| \\
 & {= 64\,{z_{{1}}}^{3}z_{{5}}-32\,{z_{{1}}}^{2}z_{{2}}z_{{4}}-16\,{z_{{1}}}^{2}{z_{{3}}}^{2}+24\,z_{{1}}{z_{{2}}}^{2}z_{{3}}-5\,{z_{{2}}}^{4}},
\end{align*}
   {which} is a basic relative invariant corresponding to the character $\nu_{  {(-3,-1)}}$,
   {while} $Q_2(z)=z_1$ is a basic relative invariant corresponding to the character $\nu_{  {(-1,0)}}$.
Therefore, we have 
$$\mathcal{ML}(P_1)=\frac{1}{(2^6 \cdot 3)^3 }  \frac{Q_1^3 }{Q_2^8}=\frac{1}{(2^6 \cdot 3)^3 }  \frac { \left( 64\,{z_{{1}}}^{3}z_{{5}}-32\,{z_{{1}}}^{2}z_{{2}}z_{{4}}-16\,{z_{{1}}}^{2}{z_{{3}}}^{2}+24\,z_{{1}}{z_{{2}}}^{2}z_{{3}}-5\,{z_{{2}}}^{4} \right) ^3 }{ z_1^8 }, $$

$$\mathcal{ML}(Q_1)=\frac{1}{(2^6 \cdot 3)^3 }  \frac{P_1^3 }{P_2^8}=\frac{1}{(2^6 \cdot 3)^3 }  \frac{(-y_{{1}}{y_{{5}}}^{3}+2\,y_{{2}}{y_{{5}}}^{2}y_{{4}}+{y_{{3}}}^{2}{y_{
{5}}}^{2}-3\,y_{{3}}{y_{{4}}}^{2}y_{{5}}+{y_{{4}}}^{4} )^3 }{y_5^8}. $$
\end{example}

\begin{rem}
Here we remark that, in general, 



\begin{equation} \label{ml of subhankel1}
 \mathcal{ML}(P_1)= \frac{1}{(2^{\ell_2 (r)+r-1} \cdot (r-1) )^{r-1} } Q_1 ^{r-1} Q_2^{-r^2+2r},
\end{equation}

\begin{equation} \label{ml of subhankel2}
 \mathcal{ML}(Q_1)= \frac{1}{(2^{\ell_2 (r)+r-1} \cdot (r-1) )^{r-1} } P_1 ^{r-1} P_2^{-r^2+2r}.
\end{equation}
%
%
\end{rem}

\section{  {Conjectures about} $b$-functions of $P_1$ and polarlizations.}
{In view of the previous section}, we   {give the} following conjecture:

\vspace{0.3cm}

{\bf Conjecture A} 
{\it 
For $P_1, P_2 , Q_1, Q_2$ defined in the   {previous} section, put
 $\tilde{Q_1} :=\displaystyle{\frac{(-1)^{r-1} }{2^{\ell_2 (r) +r-1} (r-1)^{r-1} } }Q_1$.
and
 $  {K := \phi_{(r-1, 1)}} = P_1^{r-1} P_2^{-r^2+2r}$.
Then we have 
\begin{equation} \label{b-fct of subhankel}
 \tilde{Q_1} ( \partial ) (  {K}^{s+1})= \prod_{k=1}^{r-1} (s+ \frac{k}{r-1} ) (s+ \frac{r+1}{2} )   {K}^s.
\end{equation}
{In other words}, we have $b$-functions of   {$K$} as follows:
\begin{equation} \label{b-fct of subhankel2}
b_{Q_1,   {K}} (s)  =\prod_{k=1}^{r-1} (s+ \frac{k}{r-1} ) (s+ \frac{r+1}{2} ).  
\end{equation}
}

\vspace{0.3cm}

On the other hand, in \cite{EKP} p.37,   {we find} some formulas of the Legendre transform of polarization as follows:

\begin{quote}
{\bf Example.} Let us point out an easy method of creating new functions satisfying the
projective semiclassical condition out of ones already known. It is straightforward
to compute that if $f_*$ is the multiplicative Legendre transform of $f$ on $V$ then the
multiplicative Legendre transform of the function ${\Bbb F}(x; y) = f'(x)y+f(x)$ on $V^2$ is
${\Bbb F}_* (x_*,y_*) = (d-1)^{1-d} (f'_* (y_*)x_*  - f_* (x_* ))^{d-1} 1f_* (y_* )^{2-d}$. This formula is valid also
for $d = 1$ if we agree that $0^0 = 1$.
\end{quote}
Here we call $F(x,y)$ a polarization of the homaloidal polynomial $f$. 
When $f$ is $n$-variables, $F(x,y)$ is a $2n$-variables polynomials.
{Thanks to} this result in \cite{EKP}, we can make a sequence   {of} $b$-functions.
{Namely,} starting from the sub-Hankel determinant $P_1$,
we repeat the same operations   {as above} and  we can make a sequence $\{ F_k \} _{k \geq 1}$
 of the homaloidal polynomials $F_k $ which are $k$-times polarizations 
 and their multiplicative Legendre transforms. 
Based on calcuation of examples, we pose the following conjecture:

\vspace{0.3cm}

{\bf Conjecture B}
{\it 
For   $k$-times polarization $F^{(r)}_k$ of size $r$ sub-Hankel determinant $P_1$ and its Multiplicative Legendre transform $H_k$, we have 

$$ b_{F_k, H_k} (s)=\prod_{i=1}^{r-1}  (s+\frac{i}{r-1}) (s+(r+1)2^{k-1} ). $$
}

\vspace{0.3cm}

Furthermore, we observe relative invariants of other prehomogeneous vector spaces and pose the following:

\vspace{0.3cm}

{\bf Conjecture C}
{\it 
If $F$ is an $N$-variable homaloidal polynomial of  degree $r$ and put $\mathcal{ML} (F)=H $, then the pair $(F, H)$ have a $b$-function $b^{[0]} (s)$ . Furthermore, let $F^{[k]}$ be $k$-times polarization of $F$ and its multiplicative Legendre transform $\mathcal{ML} (F^{[[k] } ) =H^{[k]} $,  the pair $(F^{[k]}, H^{[k]})$ have a $b$-function $b^{[k]} (s)$ and its form is the following:

$$b^{[k]}_{(F^{[k]}, H^{[k]})}  (s) =\prod_{i=1}^{r-1} (s+\frac{i}{r-1}) (s+N \cdot 2^{k-1}). $$

}
 \vspace{0.3cm}

Recently Conjectures B and C are claimed to be proved by F. Sato. 
The proof will be published elsewhere.

 \vspace{0.3cm}

{\bf Question 1}
{\it 
The form of $b$-function of sub-Hankel determinants are similar to the form of $b^{[k]}_{(F^{[k]}, H^{[k]})}  (s)$. Is there any relation between the spaces of sub-Hankel determinants and the spaces of the polarizations of the general homaloidal polynomials.
}
 \vspace{0.3cm}

\section{Orthogonal Polynomials and sub-Hankel determinants}

In this section, we discuss
 relations between various orthogonal polynomials and sub-Hankel determinants.  
We know a lot of relations between Hankel determinants and orthogonal polynomials, or polynomials that appeared in number theory, for example, Bernoulli polynomials, Euler polynomials, etc. 
However, for the case of sub-Hankel determinant, 
 few formulas have been recognized. 
Here we consider the following generalized Fibonacci polynomials and generalized Lucas polynomials:

${\it GFib}_n(s,t) $ is a {\it generalized Fibonacci Polynomial} defined by the recurrence 
${\it GFib}_n (s,t)=s {\it GFib}_{n-1}(s,t) +t {\it GFib}_{n-2} (s,t) $ with initial values ${\it GFib}_0 (s,t)=0$, ${\it GFib}_1(s,t)=1$.

${\it GLuc}_n(s,t) $ is a {\it generalized Lucas Polynomial} defined by the recurrence 
${\it GLuc}_n (s,t)=s {\it GLuc}_{n-1}(s,t) +t {\it GLuc}_{n-2}  (s,t)$ with initial values ${\it GLuc}_0 (s,t)=2$, ${\it GLuc}_1(s,t)=s$.

${\it NGLuc}_n(s,t) $ is a {\it normalized generalized Lucas Polynomial} defined by the recurrence 
${\it NGLuc}_n (s,t)=s {\it NGLuc}_{n-1}(s,t) +2t {\it NGLuc}_{n-2}  (s,t)$ with initial valuers ${\it NGLuc}_0 (s,t)=1$, ${\it NGLuc}_1(s,t)=s$.

In what follows, we denote by $H(x)$ the Hankel determinant $\det (x_{i+j})_{0 \leq i,j \leq r-1}$ of size $r$,
 and by $SH(x)$ the sub-Hankel determinant $\det (\hat{x}_{i+j})_{0 \leq i,j \leq r-1}$ of size $r$, where
$
 \hat{x}_{i+j} := \begin{cases} x_{i+j}, & (i+j \le r+2) \\ 0 & (i+j > r+2). \end{cases}
$
Then we have the following theorem:

\begin{thm}
{\rm (i)} 
We have $H({\it GFib}_{n+i,j}(s,t))=0$ and
\begin{equation} \label{genfib}
SH({\it GFib}_{n+i;j}(s,t))=(-1)^{ \frac{1}{2} r (r+1) } (-t)^n {\it GFib}_{n+r+1}^{r-2} (s,t). 
\end{equation}
{\rm (ii)}
We have $H({\it GLuc}_{n+i,j}(s,t))=0$ and
\begin{equation} \label{genluc}
SH({\it GLuc}_{n+i;j}(s,t))=(-1)^{ \frac{1}{2} r (r+1) +1} (s^2 +4t) {\it GLuc}_{n+r+1}^{r-2} (s,t). 
\end{equation}
{\rm (iii)} 
We have $H({\it NGLuc}_{n+i,j}(s,t))=0$ and
\begin{equation} \label{ngenLuc}
SH({\it NGLuc}_{n+i;j}(s,t))=(-1)^{ n+1+\frac{1}{2} r (r+1) } t^{n+1} 2^n (s-t+2) {\it NGLuc}_{n+r+1}^{r-2} (s,t).
\end{equation}
\end{thm}
We can prove this theorem by elementary matrix calculation. We shall write the proof in detail elsewhere.

Chebyshev polynomials of first, second, and third kind are denoted by
 $T_n (x), U_n (x)$ and $V_n (x)$ respectively. 
Note that $2T_n (x)= {\it GLuc}_n (2x, -1), U_n (x)={\it GFib}_n (2x,-1)$. 
Furthermore, Fibonacci Polynomial $F_n (x)$, Lucas Polynomial $L_n (x)$, the polynomial sequences $ \frac{x^n -y^n}{x-y} $ and $x^n +y^n$
 are also related to the generalized Fibonacci or generalized Lucas polynomials by
 $F_n (x)={\it GFib}_n (x,1), ~L_n (x)={\it GLuc}_n (x,1)$, $ \frac{x^n -y^n}{x-y} ={\it GFib}_n (x+y, -xy)$ and $x^n+y^n={\it GLuc}_n (x+y, -xy)$.
From Theorem 8, we have the following Corollary.

\begin{cor}
{\rm (i)} For Chebyshev polynomials of first kind $T_n (x)$ defined by the recurrence 
$T_n (x)=2x T_{n-1} (x) -T_{n-2} (x)$ with initial values $T_0 (x)=1, T_1 (x)=x$,  we have 
\begin{equation} \label{cheb1}
SH(T_{n+i+j} (x))=(1-x^2) T_{n+r+1}^{r-2} (x). 
\end{equation}
{\rm (ii)} For Chebyshev polynomials of second kind $U_n (x)$ defined by the recurrence 
$U_n (x)=2x U_{n-1} (x) -U_{n-2} (x)$ with initial values $U_0 (x)=1, U_1 (x)=2x$,  we have 
\begin{equation} \label{cheb2}
SH(U_{n+i+j} (x))=U_{n+r+1}^{r-2} (x). 
\end{equation}
{\rm (iii)} For Chebyshev polynomials of third kind $V_n (x)$ defined by the recurrence 
$V_n (x)=2x V_{n-1} (x) -V_{n-2} (x)$ with initial values $V_0 (x)=1, V_1 (x)=2x-1$,  we have 
\begin{equation} \label{cheb3}
SH(V_{n+i+j} (x))=-2(x-1)V_{n+r+1}^{r-2} (x). 
\end{equation}
{\rm (iv)} For Fibonacci polynomials  $F_n (x)$ defined by the recurrence 
$F_n (x)=x F_{n-1} (x) +F_{n-2} (x)$ with initial values $F_0 (x)=1, F_1 (x)=1$,  we have 
\begin{equation} \label{fibp}
SH(F_{n+i+j} (x))=(-1)^{n+\frac{1}{2} r (r+1) +1} F_{n+r+1}^{r-2} (x). 
\end{equation}
{\rm (v)} For Lucas polynomials  $L_n (x)$ defined by the recurrence 
$L_n (x)=x L_{n-1} (x) +L_{n-2} (x)$ with initial values $L_0 (x)=2, L_1 (x)=x$,  we have 
\begin{equation} \label{lunp}
SH(L_{n+i+j} (x))=(-1)^{n+\frac{1}{2} r (r+1) +1} (x^2+4)L_{n+r+1}^{r-2} (x). 
\end{equation}
{\rm (vi)} We have
\begin{equation} \label{poly1}
SH(\frac{x^{n+i+j} -y^{n+i+j} }{x-y})=(-1)^{\frac{1}{2}r (r+1)} (xy)^n \left( 
\frac{x^{n+r+1} -y^{n+r+1} }{x-y} \right)^{r-2}. 
\end{equation}
{\rm (vii)} We have 
\begin{equation} \label{poly2}
 SH( x^{n+i+j} +y^{n+i+j} )=(-1)^{\frac{1}{2} r (r+1)+1} (xy)^n (x-y)^2 (x^{n+r+1} y^{n+r+1} )^{r-2}. 
\end{equation}
\end{cor}

\begin{rem}
In the case of (iii), $V_n (x)$ does not come from generalized Fibonacci polynomials nor generalized Lucas polynomials, but, we can prove the equation (24) by the same way as other cases.
\end{rem}

{\bf Question 2}

{\it From a view point of theory of prehomogeneous vector space, the dual space has the same property as the original space. 
So we can expect for the interesting relations between basic relative invariants of dual space 
$( \mathfrak{g}, d \rho^* ,  \mathcal{SH} (r)^* )$ and orthogonal polynomials.
Can we make an interesting formula 
as above for $Q_1,Q_2$?}



\end{document}